\documentclass[11pt]{article}
\usepackage{amsfonts,latexsym}
\newcommand{\cleqn}{\setcounter{equation}{0}}
\newcommand{\clth}{\setcounter{theorem}{0}}
\newcommand {\sectionnew}[1]{\section{#1}\cleqn\clth}
\newcommand{\beq}{\begin{equation}}
\newcommand{\eeq}{\end{equation}}
\newcommand{\beqa}{\begin{eqnarray}}
\newcommand{\eeqa}{\end{eqnarray}}
\newcommand{\beaa}{\begin{eqnarray*}}
\newcommand{\eaa}{\end{eqnarray*}}
\newcommand{\nn}{\hfill\nonumber}
\newcommand{\text}{\textrm}
\newcommand \nc {\newcommand}
\nc \proof {{\em{Proof.\/}} }
\nc \qed {$\Box$\hfill}
\newtheorem{theorem}{Theorem}[section]
\newtheorem{lemma}[theorem]{Lemma}
\newtheorem{definition-theorem}[theorem]{Definition-Theorem}
\newtheorem{proposition}[theorem]{Proposition}
\newtheorem{corollary}[theorem]{Corollary}
\newtheorem{definition}[theorem]{Definition}
\newtheorem{example}[theorem]{Example}
\newtheorem{remark}[theorem]{Remark}
\nc \bth[1] { \begin{theorem}\label{t#1} }
\nc \ble[1] { \begin{lemma}\label{l#1} }
\nc \bdeth[1] { \begin{definition-theorem}\label{dt#1} }
\nc \bpr[1] { \begin{proposition}\label{p#1} }
\nc \bco[1] { \begin{corollary}\label{c#1} }
\nc \bde[1] { \begin{definition}\label{d#1}\rm }
\nc \bex[1] { \begin{example}\label{e#1}\rm }
\nc \bre[1] { \begin{remark}\label{r#1}\rm }
\nc \bcon[1] { \medskip\noindent{\it{Conjecture #1}} }
\nc \bqu[1]  { \medskip\noindent{\it{Question #1}} }
\nc {\eth} { \end{theorem} }
\nc {\ele} { \end{lemma} }
\nc {\edeth}{ \end{definition-theorem} }
\nc {\epr} { \end{proposition} }
\nc {\eco} { \end{corollary} }
\nc {\ede} { \end{definition} }
\nc {\eex} { \end{example} }
\nc {\ere} { \end{remark} }
\nc {\econ} {\smallskip}
\nc {\equ} {\smallskip}
\nc \eqref[1] {{\rm{(\ref{#1})}}}
\nc \thref[1]{Theorem \ref{t#1}}
\nc \leref[1]{Lemma \ref{l#1}}
\nc \prref[1]{Proposition \ref{p#1}}
\nc \coref[1]{Corollary \ref{c#1}}
\nc \deref[1]{Definition \ref{d#1}}
\nc \exref[1]{Example \ref{e#1}}
\nc \reref[1]{Remark \ref{r#1}}
\def \CM  {C^\infty (M)}
\def \CB  {C^\infty (B)}
\def \pcb  {\pi^*(C^\infty (B))}
\def \h  {\hbar}
\def \Wh  {W_{\hbar}}
\def \WhD  {W^{\hbar}_D}

\def \WB {W_B}
\def \Wp {W_{\pi}}
\def \WM {W_M}
\def \WBh  {W_{B, \hbar}}
\def \Wph  {W_{\pi, \hbar}}
\def \WMh  {W_{M,\hbar}}
\def \W  {W}
\def \p {\partial}
\def \a {\alpha}
\def \b {\beta}
\def \om {\omega}
\def \sig {\sigma}
\def \ro {\rho}
\def \Ga {\Gamma}
\def \ga {\gamma}
\def \La {\Lambda}
\def \de {\delta}
\def \Om {\Omega}
\def \ph {\varphi}

\def \g {\gamma}
\def \i {i}
\def \bI {\bar{I}}
\def \bp {\bar{\varphi}}
\def \t {{\mathcal T}}
\def \O {{\mathcal O}}
\def \O {{\mathcal O}}
\def \U {{\mathcal U}}
\def \J {{\mathcal J}^\infty}

\def \Rset {{\mathbb R}}

\def \Zset {{\mathbb Z}}

\def \ra {\rightarrow}
\def \wt {\widetilde}
\def \ol {\overline}

\def \hra {\hookrightarrow}
\def \st {\ast}
\def \stw {\bullet}
\def \pe {\star}
\def \pew {\bullet}
\def \pol { {\mathrm{pol}} }
\def \fib { {{\mathrm{fib}}} }
\def \id { {\mathrm{id}} }
\def \ad { {\mathrm{ad}} }

\def \span { {\mathrm{span}} }
\def \const { {\mathrm{const}} }

\def \mod { {\mathrm{mod}} }

\renewcommand \ker { {\mathrm{Ker}} }
\begin{document}
%%%%%%%%%%%%%%%%%%%%%%%%%%%%%%%%%%%%%%%%%%%%%%%%%%%%%%%%%%%%%%%%%%%%%%%
%%%%%%%%%%%%%%%%%%%%%%    Title    %%%%%%%%%%%%%%%%%%%%%%%%%%%%%%%%%%%%%%%%%%%
\title{{\LARGE\bf{
Deformation Quantization \\ of Lagrangian Fiber Bundles
}}}
\author{
Nicolai~Reshetikhin and Milen~Yakimov \\
{\normalsize{Department of Mathematics,}} \\
{\normalsize{University of California at Berkeley,}} \\
{\normalsize{Berkeley, CA 94720, USA}}
\thanks{E-mail addresses: reshetik@math.berkeley.edu, 
                          yakimov@math.berkeley.edu}
\\ \hfill \\
{\it{Dedicated to the memory of Mosh\'e Flato}} 
}
\date{}
\maketitle
\begin{abstract}
Let $(M, \om)$ be a symplectic manifold. A Lagrangian fiber bundle 
$\pi: \: M \ra B$ determines a completely integrable system on $M.$ 
First integrals of this system are the pull--backs
of functions on the base of the bundle.  We show that for each  
Lagrangian fiber bundle $\pi$  there exist star products on $\CM[[\h]]$
which do not deform the pointwise multiplication on the subalgebra
$\pcb[[\h]].$ The set of equivalence classes of such star products is in
bijection with formal deformations of the symplectic structure $\om$
for which $\pi: M \ra B$ remains Lagrangian taken modulo formal
symplectomorphisms of $M.$
\end{abstract}

%%%%%%%%%%%%%%%%%%%%   Introduction   %%%%%%%%%%%%%%%%%%%%%%%%%%%%%%%%%%%%%%%%
\sectionnew{Introduction}
%%%%%%%%%%%%%%%%%%%%%%%%%%%%%%%%%%%%%%%%%%%%%%%%%%%%%%%%%%%%%%%%%%%%%%%%%%%%%%
%First, let us explain the motivation of the problem.
The phase space of a classical Hamiltonian system is a symplectic 
manifold $(M, \om).$ Smooth functions on the phase space form classical
algebra of observables. It is a Poisson algebra with Lie bracket given
by
\beq
\{f, g \} = \om^{-1}(d f \wedge d g), \quad f, g \in \CM
\label{0.lie}
\eeq
$(\om^{-1} \in \wedge^2 TM$ is the bivector field dual to the two-form
$\om)$ and commutative product the pointwise product on $\CM.$

A completely integrable classical Hamiltonian system on a 
$2 n$-dimensional symplectic manifold admits $n$ independent integrals
which commute with respect to the Poisson bracket. Geometrically this means
that away from the singular level surfaces the symplectic manifold $M$ is
fibered over an open subset of  $\Rset^n.$  More generally for any fiber
bundle $\pi: M \ra B,$  whose fibers are Lagrangian submanifolds of 
$(M, \om)$ the subalgebra $\pcb$ of $\CM$ is Poisson commutative. This a
geometric model for a completely integrable classical system without
singularities.

Quantization as a way to pass from classical Hamiltonian mechanics to
quantum mechanics mathematically can be regarded as a procedure of
replacing the classical algebra of observables (the Poisson algebra of
functions on a symplectic manifold) by an associative noncommutative
algebra of quantum observables. One approach to quantization is to
consider the algebra of quantum observables as a (formal) deformation
of the Poisson algebra $(\CM, \{., .\})$ of classical observables in the
category of associative algebras. For precise definitions see Sect.~2.
This approach is known as  deformation quantization. It was pioneered by
H.~Weyl \cite{We}, and J.~Moyal \cite{M}, who considered the case of
symplectic vector spaces and developed further by J.~Vey \cite{V},
F.~Berezin
\cite{B} and F.~Bayen, M.~Flato, C.~Fronsdal, A.~Lichnerowicz, and
D.~Sternheimer \cite{BFFLS} in the case of general Poisson manifolds.
There were several techniques employed in the study of the deformations
of the Poisson algebra $(\CM, \{., .\}).$
In \cite{V, BFFLS} a cohomological technique, based on results of \cite{G}
was developed. It was successively used in \cite{DL} to prove that any
symplectic manifold can be quantized. Fedosov \cite{F} found a geometric
construction of deformations in the symplectic case.
Finally the problem of quantization was completely solved by Kontsevich
\cite{K} for any Poisson manifold. He derived it from the more general
result of formality of the differential graded Lie algebra of the
Hochschild complex of functions on any smooth manifold. 

For a Lagrangian fiber bundle $\pi : M \ra B$ inside the algebra of
classical observables $\CM$ sits the Poisson commutative subalgebra of
classical integrals of motion $\pcb.$ It is natural to ask two questions
in this situation:
\begin{enumerate}
\item Is it possible to construct a quantization of the Poisson algebra 
$(\CM, \{., .\})$ such that the subalgebra $\pcb$ will not deform
and thus will form an algebra of quantum integrals of motion?
\item If it is possible, then ``how many of such quantizations'' one
can construct?
\end{enumerate}

In this paper we give answers to both questions using Fedosov's technique. 
Our results are  summarized in the following Theorem.

\bth{Summ} {\em{(}}i{\em{\/)}} The classes of star products on $(M, \om)$
which are equivalent to products keeping the commutative algebra $\pcb$
undeformed are in one to one correspondence with the deformations of the
two-form $\om$ in $H^2(M)[[\h]]$ having  representatives in $Z^2(M)[[\h]]$
with respect to which the bundle $\pi : M \ra B$ remains Lagrangian.

{\em{(}}ii{\em{\/)}} Any star product on $(M, \om)$ for which the space
$\pcb[[\h]]$ generates a commutative subalgebra of $\CM[[\h]]$ is
equivalent to one that keeps the usual product on $\pcb[[\h]]$ undeformed.
\eth

This Theorem provides a positive answer to Question 1. In the framework
of Fedosov's quantization the formal two-form associated to a class
of star products on $(M, \om)$ with the decsribed properties
is its characteristic form (class).
In view of \cite{Moser}, the discussed classes of star products on
$(M, \om)$ are in one to one correspondence with formal
$(\Rset[[\h]])$ deformations of the symplectic structure $\om$ on $M$
for which $\pi : M \ra B$ remains Lagrangian, taken modulo the action of
the group of formal symplectomorphisms of $(M, \om)$ starting with the
identity diffeomorphism of $M.$

\bre{QuantCorr} Often a quantization $\pe$ of $(M, \om)$ is apriori
given (e.g. quantization of a cotangent bundle in terms of differential
operators). In this situation one wants to deform
(if it is possible) the Poisson commutative algebra
$\pcb$ to a commutative subalgebra of $(CM[[\h]], \pe)$
by introducing {\em{quantum corrections}} to the commuting quantities.
For this problem \thref{Summ} can be used in the following manner.
The first step is to check whether the characteristic form of the product
$\pe$ is of the required type. If it is so, let $\st$ be an equivalent
to $\pe$ product which keeps the algebra $\pcb$ undeformed and
the corresponding formal equivalence be given by the differential
operator $P = \id + O(\h).$ (An iterative procedure for constructing
such $\st$ is contained in Sect.~5.)
Then $(P^{-1}(\pcb)[[\h]], \pe)$ is a commutative subalgebra of
$(\CM[[\h]], \pe)$ and the quantum correction to each classical
integral of motion $f \in \pcb$ is $(P^{-1} -\id)(f).$
\ere

The paper is organized as follows. In Sect.~2 we remind some basic
definitions and facts about Fedosov's quantization. Sect.~3 contains the
construction of special connections which we use in Sections~4 and~5. The
proof of the existence of associative
star products on $\CM[[\h]]$ preserving given Lagrangian fibration is
presented in Sect.~4. Sect.~5 is devoted to the classification of such star
products. \thref{Summ} is a combination of \thref{3.1} and \thref{4.1}. 

One can restrict to the set of star products on $(M, \om)$ that
keep $\pcb[[\h]]$ undeformed and consider a finer equivalence
relation given by differential operators that preserve
$\pcb[[\h]].$ A classification result in this setting
as well as extensions to the case of noncommutative integrable
systems will be discussed in a separate publication.

%In the Conclusion we discuss possible directions in quantization of
%integrable systems.

{\flushleft{\bf{Acknowledgements}}}

\medskip\noindent
We benefited from discussions with B.~Fedosov, A.~Givental, 
L.~Takhtajan, and  A.~Weinstein whom we thank for this. 
We thank also Technishe Universit\"at Berlin for the warm hospitality
when this work was being completed. The research of N.R. was partially
supported by NSF grant DMS96-03239,
and M.Y. received support from NSF grants DMS96-03239 and DMS94-00097.   
%%%%%%%%%%%%%%%%%%%%%%%%%%%%%%%%%%%%%%%%%%%%%%%%%%%%%%%%%%%%%%%%%%%%%%%%%%%%%%%
\sectionnew{Deformation Quantization and Fedosov's Construction}
%%%%%%%%%%%%%%%%%%%%%%%%%%%%%%%%%%%%%%%%%%%%%%%%%%%%%%%%%%%%%%%%%%%%%%%%%%%%%%%
\subsection{Basic Notions}
%First we briefly recall the main definitions for Deformation
%Quantization. For more details we refer to \cite{BFFLS, W}. In
%Subsect. 2.2 and 2.3 we recollect some major facts on the classification
%of star products on symplectic manifolds based on Fedosov's work
%\cite{F}.
%They will be used in Sections 4 and 5, which contain the main results of
%this article.   

\bde{1.star}A formal local star product on a Poisson manifold 
with Poisson bracket $\{., .\}$ is an associative
$\Rset[[\h]]$-linear product on $\CM[[\h]]$ of the type
\beq
f \st g = fg +\sum_{l=1}^\infty \h^l Q_l(f,g), 
\quad f, \; g \in \CM[[\h]],
\label{1.st}
\eeq
where 
\begin{enumerate}
\item $Q_l$ are local (i.e. bidifferential) operators,
\item $f \st g - g \st f = \h \{f, g\} + O(\h^2),$ 
$\forall f, g \in \CM[[\h]],$
\item $1 \st f= f \st 1=f,$ i.e $Q_l(1, f)= Q_l(f, 1)=0,$
$\forall f \in \CM[[\h]],$ $l \geq 1.$
\end{enumerate}
\ede

In this paper we will consider only {\em{formal local}} star products
and we will call them just star products.

\bde{1.s_equ}Two star products $\st_1$ and $\st_2$ are called equivalent
if there exists a local (i.e. differential) $\Rset[[\h]]$-linear
operator $P: \CM[[\h]] \ra \CM[[\h]]$ of the type
\beq
P = \id + \sum_{l=1}^\infty \h^l P_l
\label{1.s-equ}
\eeq
for some differential operators $P_l$ on $M$ such that
\beq
P(f\st_1 g) = P f \st_2 P g, \quad 
\forall f, g \in \CM.
\label{1.equiv}
\eeq
\ede

For more details on these definitions we refer to \cite{BFFLS, M}.
%%%%%%%%%%%%%%%%%%%%%%%%%%%%%%%%%%%%%%%%%%%%%%%%%%%%%%%%%%%%%%%%%%%%%%%%%
\subsection{The Weyl Bundle}
%This Subsection contains an overview of
%Fedosov's
%construction of star products on symplectic manifolds. 

Let  $(M, \om)$ be a symplectic manifold of dimension $2 n.$ The Weyl
bundle $W$ on $M$ is defined by
\beq
\W=S(T^* M), 
\label{1.1}
\eeq
where $S(V)$ denotes the completed symmetric power 
$S(V)= \coprod_{l=0}^\infty S^l(V)$ of a vector space $V.$
For any point $x \in M$ the commutative algebra $S(T_x^* M)$ can be
naturally identified with the algebra $\J_0(T_x M)$ of $\infty$-jets of
real valued functions on $T_x M$ at 0.
%The algebra of global sections of $W$ equipped with the standard
%multiplication in $S(V)$ is naturally identified with the
%algebra of functions on a formal neighborhood of the 0 section of $T M.$
More often we will use the bundle
\beq
\Wh=S(T^* M)[[\h]]=\W[[\h]],
\label{1.1a}
\eeq
where $\h$ is a formal variable and will still call it Weyl bundle 
of $M$.

The tangent space $T_x M$ at any $x \in M$ is a symplectic vector space
and the Moyal-Weyl product equips any fiber $(\Wh)_x$ of $\Wh$
with a structure of noncommutative associative algebra. It
becomes a graded algebra if we put
\beqa
  &&\deg v = l, \quad \mbox{for } v \in S^l(T_x^*M)
    \label{1.grad}\\
  &&\deg \h =2.
    \label{1.gradh}
\eeqa
The homogeneous component of $(\Wh)_x$ in degree $l$ will be denoted by
$(\Wh^l)_x.$ The fiberwise product on $\Wh$ determines an
associative product on the space $\Ga(\Wh)$ of global sections of
$\Wh,$ denoted by $\circ.$

In local coordinates $(x^1, \dots, x^{2n}): U \ra \Rset^{2n}$ the
symplectic
form $\om$ can be written as 
\[
\om = \om_{jl} \; d x^j \wedge d x^l.
\]
Here and later we use the standard convention of summation
over repeated indices.
Denote by $y^j$ the image of $d x^j$ in $S(T_x^*M).$ Any section
$a\in
\Ga(\Wh, U)$ is a formal power series:
\beq
a= \sum_{l\geq 0} \sum_\a a_{l, \a}(x) y^\a \h^l,
\label{1.loc}
\eeq
where $\a=(\a_1, \ldots, \a_{2n})$ is a a multiindex. In this notation
the Moyal--Weyl product of two sections $a, \; b \in \Ga(\Wh, U)$ is

\beq
a \circ b = \exp \left( \frac{\h}{2} \; \om^{jl}(x) \; 
                        \frac{\p}{\p y^j}
                        \frac{\p}{\p z^l}  \right)
                        a(y) b(z)\Big\vert_{z=y}.
\label{1.3}
\eeq
\subsection{Constructing Fedosov's connections}
The main idea of Fedosov's quantization is to construct a flat connection
in the Weyl bundle $\Wh$ for which the exponential map identifies 
flat sections with smooth functions on $M.$ Then the Moyal--Weyl product
\eqref{1.3} descends to a star product on $\CM[[\h]].$

Denote by $\La$ the exterior algebra bundle on $M:$
\beq
\La = \wedge(T^*M).
\label{1.ext}
\eeq
Moyal--Weyl product on $\Wh$ together with the exterior product
on $\La$ equips $\Wh\otimes\La$ with a structure of a super algebra.
($\Wh$ sits in even degree, $\La^l$ as usual in even/odd degree
depending on the parity of $l$.) The associative products in 
$(\Wh\otimes\La)_x,$ $x \in M$ and $\Ga(\Wh\otimes\La)$ will
be again denoted by $\circ,$ and the corresponding supper bracket
by $[.,.].$

There are two canonical fiberwise endomorphisms of $\Wh\otimes\La,$
defined locally by
\beq
   \de(a) =\sum_j d x^j \wedge \frac{\p a}{\p y^j}
\label{1.de}
\eeq
and
\beq
\de^{-1} (a) = 
  \left\{
  \begin{array}{cc}
      \frac{1}{p+q} \sum_j \: y^j \left(\frac{\p}{\p x^j} 
%                                  \rfloor
                                   \lrcorner  \: a \right),
      & \mbox{for $a \in (\Wh^p\otimes \La^q)_x, \; p+q>0,$} \\
      0,
      & \mbox{for $a \in (\Wh^0\otimes \La^0)_x.$}
  \end{array}
  \right.
\label{1.de-1}
\eeq
They satisfy
\[
\de^2 = (\de^{-1})^2 =0
\]
and a Hodge type relation 
\[
a = (\de \de^{-1} + \de^{-1} \de)a + a_{00}, \quad
\forall a \in (\Wh \otimes \La)_x,
\]
where $a_{00}$ is the component of $a$ in
$(\Wh^0 \otimes \La^0)_x.$

Fix a symplectic torsion free connection on $M.$ Let $\p$ denote its
covariant derivative. It extends in a natural way to a derivation of
$(\Ga(\Wh\otimes\La), \circ).$ In local Darboux coordinates
on an open $U$
\beq
\p a = d a +\frac{1}{\h}[\Ga, a], \quad a \in \Ga(W\otimes\La, U),
\label{1.11}
\eeq
where
\beq
\Ga =\frac{1}{2} \Ga_{jkl}y^j y^k d x^l
\label{1.12}
\eeq
and $\Ga_{jkl}$ are the Christoffel symbols of the connection in the
considered coordinate system.

Consider more general derivations of $(\Ga(\Wh\otimes\La), \circ)$
(or ``nonlinear connections'' on $\Wh$) of the type
\beq
D = \p - \de + \frac{1}{\h}\ad(\g)
\label{1.13}
\eeq
with
\[
\g \in \oplus_{p \geq 3}\Ga(\Wh^p\otimes\La^1).
\]
One computes 
\beq
D^2 a = \frac{1}{\h}[\Om, a], \forall 
    a \in \Ga(\Wh\otimes\La)
\label{1.14}
\eeq
with
\beq
\Om = R +\p \ga - \de \ga + \frac{1}{\h} \ga^2 + \om
\nn
\eeq
called Weyl curvature of $D$ (the term $\om$ is for convenience only).
Here $R$ is essentially the curvature of
$\p$
\beq
R = \frac{1}{4} R_{ijkl} y^i y^j d x^k \wedge d x^l, \quad 
R_{ijkl} = \om_{im} R_{jkl}^m.
\label{1.15}
\eeq

The connection $D$ is called a Fedosov's connection when it is flat
$(D^2=0).$ That is when
\beq
\Om \in Z(\Ga(\Wh\otimes\La)),
\label{1.16}
\eeq
where $Z(\Ga(\Wh\otimes\La))$ denotes the supercenter of 
$(\Ga(\Wh\otimes\La), \circ).$ It is clear that
\beq
Z(\Ga(\Wh\otimes\La))=\Om(M)[[\h]],
\label{1.17}
\eeq
where $\Om(M)= \Ga(\La, M)$ is the space of differential forms on $M.$
So eq. \eqref{1.16} is equivalent to $\Om \in \Om^2(M)[[\h]].$ In that
case
Bianchi identities together with eq. \eqref{1.14} imply
that $D \Om = d \Om =0,$ i.e.
\[
\Om \in Z^2(M)[[\h]],
\]
where now $Z(M)$ stays for closed forms on $M.$

\bth{1.1} {\em{(}}Fedosov{\em{\/)}} For any choice of the Weyl curvature 
$\Om \in Z^2(M)[[\h]],$ $\Om =\om +O(\h),$ there exists
a unique Fedosov's connection on $W$ with \\ 
$\ga \in \oplus_{p \geq 3}\Ga(\Wh^p\otimes\La^1)[[\h]]$ and
$\de \ga = 0.$ The components of of the one-form $\g$ in the grading 
\eqref{1.grad}--\eqref{1.gradh} of $\Wh$ are consecutively determined from
\beq
\g = \de^{-1}(-\Om + \om +R +\p \g + \frac{1}{\h}\g^2).
\label{1.18}
\eeq
\eth

\subsection{Construction of Fedosov's star products}
For a Fedosov's connection $D$ on $\Wh$ the space 
\beq
\WhD = \{ a \in \Ga(\Wh) | \; D a =0 \}
\label{1.10}
\eeq
of flat sections of $\Ga(\Wh)$ is a subalgebra of $(\Ga(\Wh), \circ).$ 

Given a section $a \in \Ga(\Wh)$ denote by $a_0$ its component
in $S^0(T^* M)[[\h]] = \CM[[\h]],$ i.e. the constant (in $y$)
term in \eqref{1.loc}.

\bth{1.2}{\em{(}}Fedosov{\em{\/)}} For any $f \in \CM[[\h]]$
there exists a unique $\sig(f) \in \WhD$ such that
$\sig(f)_0 = f.$ The graded components $\sig(f)$ with respect to the
grading \eqref{1.grad}--\eqref{1.gradh} of
$\Wh$ are determined inductively by
\beq
\sig(f) = f + \de^{-1}\left( \p(\sig(f)) 
+\frac{1}{\h}[\g, \sig(f)]\right).
\label{1.20}
\eeq
\eth

The above theorem allows to introduce an associative product
on $\CM[[\h]]$ by
\beq
f \st g = \left(\sig(f) \circ \sig(g) \right)_0,
\quad f, g \in \CM[[\h]].
\label{1.21}
\eeq
One easely sees that it is a deformation of the Poisson structure
on $(M, \om).$
%%%%%%%%%%%%%%%%%%%%%%%%%%%%%%%%%%%%%%%%%%%%%%%%%%%%%%%%%%%%%%%%%
\subsection{Classification of star products}
The classification of star products on a symplectic 
manifold was obtained in \cite{DL} via cohomological methods developed
in \cite{BFFLS}. In Fedosov's approach it was done in \cite{F2, NT, X}.
For comparison between the two approaches see \cite{D}. 
The classification is summarized in the following theorem.

\bth{1.3} {\em{(}}i{\em{\/)}} There exists a bijective correspondence
between the set of isomorphism classes of star products on $(M, \om)$ and
the classes of formal deformations of the two-form $\om$ on
$M$ $(\equiv \{ [\Om] \in H^2(M)[[\h]] | \Om =\om + O(\h)\}).$

{\em{(}}ii{\em{\/)}} Any star product is equivalent to a Fedosov's one
and the class of this product referred to in part {\em{(}}i{\em{\/)}} is
the
class of the curvature form $[\Om]$ of the corresponding Fedosov's
connection {\em{(}}see \thref{1.1}{\em{\/)}}.
\eth

Using Moser's classification of nearby symplectic structures 
\cite{Moser}, \thref{1.3} can be rephrased as follows:

{\em{The equivalence classes of star products on $(M, \om)$ are in one to
one correspondence with the equivalence classes of formal $(\Rset[[\h]])$
deformations of the symplectic structure $\om$ on $M.$ }} \\
\noindent
The equivalence relation in the second set is defined by the action of
the group of formal symplectomorphisms of $(M, \om)$ starting with the
identity diffeomorphism of $M.$ 

\subsubsection{Semiclassical Fedosov's exponential map}
Consider the formal bundle $\J(M)$ over $M$ whose fiber over $x \in M$ is
the set of $\infty$-jets $\J_x(M)$ of real valued functions on $M$ at $x.$
For a symplectic manifold $(M, \om)$ each fiber of $\J(M)$ has a  
natural Poisson structure and so does each fiber of the
Weyl bundle
$W$. Emmrich and Weinstein
\cite{EW} found that a semiclassical analog of Fedosov's construction
gives a fiberwise isomorphism
\beq
\exp: \J(M) \ra W.
\label{1.iso}
\eeq
Their construction goes as follows. Consider a symplectic torsion free
connection $\p$ on $M.$ There exists 
\beq
\ga \in \oplus_{p \geq 3} \Ga(W^p\otimes \La^1)
\nn
\eeq
such that
\beq
D^0 = \p -\de + \{\ga, .\}_{\fib}
\label{1.22}
\eeq
is a flat connection on $W,$ i.e. $(D^0)^2=0.$
With $\{.,.\}_{\fib}$ we denote the fiberwise Poisson bracket on 
$W.$ One such $\ga$ can be inductively computed from
\beq
\ga = \de^{-1}\left(R + \p \ga +\frac{1}{2}
\{\ga, \ga \}_{\fib}\right)
\label{1.23}
\eeq
(cf. \eqref{1.18}). Here $R$ is the curvature of $\p,$ as in eq.
\eqref{1.15}.
Denote the jet of flat sections of $\W$ with respect to $D^0$ at
$x \in M$ by
\beq
\J_x(W)_{D^0}=\{ a \in \J_x(W)| \; D^0 a=0\}.
\label{1.24}
\eeq
For any $f \in \J_x(M)$ there exists a unique $\bar{f} \in \J_x(W)_{D^0}$
such that
\[
\left(\bar{f}\right)_0=f.
\]
Iteratively it is computed from 
\beq
\bar{f} = f + \de^{-1}\left( \p(\bar{f})
        + \{\ga, \bar{f}\}\right)
\label{1.sec}
\eeq
cf. \eqref{1.20}. At the end one can define 
\beq
\exp_x(f) = \bar{f}(x).
\label{1.25}
\eeq
  {}From eq. \eqref{1.sec} it follows that the evaluation in the RHS
makes sense. 

\bpr{we}
The map $\exp: \J_x(M) \ra (W)_x$ defined by \eqref{1.25} is
a fiberwise Poisson isomorphism.
\epr
\proof From the fact that $D^0$ is a differentiation of the Poisson
algebra $\Ga(W)$ and thus differentiation of the algebras $\J_x(W)$ for
any $x \in M$ it follows that $\exp$ is a fiberwise Poisson map. The map
$\exp$ is a fiberwise isomorphism because 
\beq
\exp_x(f) = f + \frac{\p f}{\p x^j} y^j + O(y^2)
\label{1.26}
\eeq   
in any coordinate system $(x^1, \ldots, x^{2 n})$ in a neighborhood
of $x \in M.$
\hfill \qed

\subsubsection{A proof of the classification Theorem} 
Here we explain some parts of the proof of \thref{1.3}. 
In particular following  P.~Xu \cite{X} we show how one obtains from
a given star product $\st$ on $(M, \om)$ a Fedosov's star product
equivalent to it. The product $\st$ is a local one and thus it induces
a star product on $\J_x(M),$ $x \in M.$ Using the Emmrich--Weinstein
exponential map $\exp$ we get a fiberwise star product $\stw$ on $\Wh:$
\beq
a \stw b := \exp_x( \exp^{-1}_x(a) \st \exp^{-1}_x(b) ),
\; \mbox{for }a, \: b \in (\W)_x, (x \in M).
\label{1.Wprod}
\eeq

This product is a deformation quantization of the fiberwise Poisson
algebra structure $\{., .\}_\fib$ on $W$ due to \prref{we}.
It is fiberwise equivalent to the Moyal--Weyl product
$\circ$ \eqref{1.3}, because $T M$ is a regular Poisson manifold with
contractable symplectic leaves. Therefore there exists a differential
operator of the form
$$ P = 1+\sum_{l=1}^\infty \h^l P_l,$$
where $P_l$ are differential operators on $T M$ that involve derivatives
along the tangent spaces only, such that
\[
 P(a\stw b)=P(a) \circ P(b), \quad
 \forall a, b \in \Ga(\Wh).
\]
Let  
\[
j: \CM \ra \Ga(\J(M))
\]
be the natural evaluation map, sending a smooth function to the
$\infty$-jet that it determines at each point of $M.$ Then the
composition 
\beq
\ro = P \cdot \exp \cdot j : (\CM[[\h]], \st) \ra (\Ga(\Wh), \circ)
\label{1.com}
\eeq
is an injective homomorphism having the following properties 

(i) $\left(\ro(f)\right)_0 = f,$

(ii) $\ro(f) \equiv f + \de^{-1} d f \; ( \mod \; \Ga (\Wh^2)).$

Any Fedosov's connection $D$ gives rise to a map $\sig_D$ defined in
\thref{1.2}, which satisfies (i)--(ii). P.~Xu proved a converse
statement.

\bth{1.5} Let $\st$ be a star product on $M$ and 
$\ro: (\CM[[\h]], \st) \ra (\Ga(\Wh), \circ)$ 
be a homomorphism that satisfies (i)--(ii).
Then there exists a Fedosov's connection $D$ such that
$$ \ro = \sig_D.$$
\eth

This Theorem applied to the composition \eqref{1.com} gives a Fedosov's
connection $D$ on $\Wh$ and a Fedosov's star product $\st_{F}$ on $M$
equivalent to the initial one $\st$. The class of the curvature form of
$D$ is called characteristic class of the star product $\st$ and is the
two-form that appears in \thref{1.3}.

Finally two Fedosov's products are equivalent if and only if their 
characteristic classes are equal. We will not need a proof of this
fact here. The reader can find it in \cite{F2, NT}.
%%%%%%%%%%%%%%%%%%%%%%%%%%%%%%%%%%%%%%%%%%%%%%%%%%%%%%%%%%%%%%%%%%%%%%%
\sectionnew{Lagrangian fiber bundles and some connections associated
with them}
%%%%%%%%%%%%%%%%%%%%%%%%%%%%%%%%%%%%%%%%%%%%%%%%%%%%%%%%%%%%%%%%%%%%%%%
Let $(M, \om)$ be a symplectic manifold and $\pi: M \ra B$ be a
Lagrangian fiber bundle, i.e $\pi$ is a fiber bundle whose fibers are
Lagrangian submanifold of $(M, \om).$ We assume also that the fibers of
$\pi$ are connected and that for each $f \in \CB$ the Hamiltonian vector
field $H(\pi^* f)$ is complete. Equivalently the fibers of $\pi$ should
be diffeomorphic to $\Rset^n/\Zset^k$ for some $k \leq n.$ 
We do not require compactness of the fibers thus $k$ can be $< n.$

%%%%%%%%%%%%%%%%%%%%%%%%%%%%%%%%%%%%%%%%%%%%%%%%%%%%%%%%%%%%%%%%%%%
\subsection{Action--angle coordinate charts}
In this subsection we review the construction of action--angle
coordinate charts on $M$ (see \cite{A, GS, Du}). 

For any $b \in B$ denote
\beq
M^b= \{x \in M|\; \pi(x) = b\}.
\label{2.1}
\eeq
$T_b^*B$ acts on $M^b$ in the following way. Let $f$ be a smooth function
defined in a neighborhood of $b.$ The action of $(d f)_b$ on $ x \in
M^b$ is defined by
\beq
\tau((df)_b)x = g_1^{H(\pi^* f)}(x),
\label{2.2}
\eeq
where on the RHS $g$ denotes the flow on $M$ corresponding to a
Hamiltonian vector field (upper index) at a given time (lower index). 
This flow preserves $M^b$ and
thus it depends on $(d f)_b$ only, so the action is correctly defined. 
$T^*_bB$ is a commutative group, thus the stabilizers of all points of
$M^b$ coincide. Denote the corresponding subgroup of $T^*_bB$ by $P_b$. It
is obvious that $P_b \cong \Zset^k$ for some $k \leq n.$ Varing $b$ we get
a submanifold $P$ of $T^*B.$ It is also clear that it is a Lagrangian
submanifold of $T^*B.$ The action of $T^*_bB$ on $M^b$ defines a faithful
action of the fibers $T^*_b B /P_b$ of $T^*B/P$ on the fibers
$M^b$ of the bundle $\pi.$ 

Let $U \subset B$ be an open set for which there exists a section
\beq
s: \; U \ra M
\label{2.s}
\eeq 
of $\pi$ $(\pi \circ s =\id),$ such that

\begin{enumerate}
\item $s(U)$ is a Lagrangian submanifold of $M$
\item  $s$ is a diffeomorphism onto its image.
\end{enumerate}
The two conditions roughly mean that we consider a Lagrangian submanifold
$s(U)$ transversal to the fibration. Any such section determines a
symplectic diffeomorphism
\beq
\pi^{-1}(U) \cong (T^*B/P)|_U,
\label{2.4}
\eeq
where $T^*B$ is equipped with the standard symplectic structure
as a cotangent bundle. To prove this simply act with $T^*_bB/P_b$ on
$s(b)$ for $b \in U.$ 

In the case when $U$ is diffeomorphic to a domain in $\Rset^n$
trivializing $(T^*B/P)|_U$ gives a chart on
$\pi^{-1}(U)$ of the type
\beq
  (I^1, I^2, \ldots, I^n,
  \ph^1, \ph^2, \ldots, \ph^n): \;
  \pi^{-1}(U) \ra O \times(\Rset^n/\Zset^k).
\label{2.5}
\eeq
Here $O$ is an open domain of $\Rset^n$, which from now on will be
denoted by the same letter $U$. The lattice $\Zset^k$ is assumed
to be embedded in the standard way in $\Rset^k \subset \Rset^n$ 
spanned by the first $k$ coordinate vectors.
The coordinates $(I, \ph)$ are called action--angle coordinates, 
although in the case $k<n$ only $\ph^1, \ldots, \ph^k$ are
angles. Each $b\in B$ has a neighborhood of the
above type and thus $M$ can be covered by action--angle charts.

Let
\beq
(\bar{I}^1, \bar{I}^2, \dots, \bar{I}^n,  
\bar{\ph}^1, \bar{\ph}^2,\ldots, \bar{\ph}^n):
\pi^{-1}(\ol{U}) \ra \ol{U} \times (\Rset^n/\Zset^k).
\label{2.6}
\eeq
be another action--angle chart on $\pi^{-1}({\ol{U}})$ for some
open
domain $\ol{U}$ with $U\cap \ol{U} \neq \emptyset.$ The change of
coordinates on the intersection of the two action--angle charts is of a
particularly simple form:
\beqa
  &&\bar{I}^\a = f^\a(I), \quad \a =1, \dots, n,
                                       \label{2.7} \\
  &&\bar{\ph}^j= A_l^j(I) \ph^l + B^j(I), 
    \quad i, \a =1, \ldots, n,
                                        \label{2.8}
\eeqa
where $A_l^j(I)$ is the inverse of the matrix 
$\frac{\p \bar{I}^l}{\p I^l}$ and $f^\a, B^1, \ldots, B^n$ are
some $C^\infty$ functions.

The RHS of \eqref{2.8} is linear in $\ph^l$ because of
the way we defined the identification \eqref{2.4}. $B^j(I),$
$j=1, \ldots, n$ are the coordinates of the Lagrangian section
$s|_{U\cap\ol{U}}$ of $\pi$ (see \eqref{2.s}) in the action--angle
chart $(\bar{I}, \bar{\ph})$ on $\pi^{-1}(U\cap\ol{U}):$
\[
s(\bar{I})=\left(\bar{I}, B^j\left(I(\bar{I})\right)\right).
\]
This imposes some standard differential conditions on the functions
$B^j(I).$ 

Note also that \eqref{2.8} should define a map from $\Rset^n/\Zset^k$
to $\Rset^n/\Zset^k.$ This means that $A^j_l(I)+B^j(I)$ should be integers
for $j=1, \ldots n,$ $l=1, \ldots k$ and thus constants as functions of
the action variables $I.$

\subsection{Special symplectic connections}

In the second part of of this section we consider certain type of
symplectic connections associated to a Lagrangian fiber bundle
$\pi: M \ra B.$ They will be used in Sections 4 and 5 as initial
connections $\p$ for special Fedosov's connections.

First we give a computational proof of our main theorem and then explain
some parts of it in geometric terms.

\bth{2.1} There exists a torsion free symplectic connection $\Ga$ on $M$
such that in each action--angle chart $(I, \ph)$ as in
\eqref{2.5} the Christoffel symbols of $\Ga$ satisfy the following
properties:
\beqa
&& \Ga_{\ph^j \ph^l}^{I^\a}=0, 
   \; \mbox{i.e. }\Ga_{\ph^j \ph^l \ph^m}=0,
\label{2.12}\\
&& \Ga_{ \ph^j \ph^l }^{ \ph^m }=
   \Ga^{I^\b}_{I^\a \ph^j }=0, \;
   \mbox{i.e. }
   \Ga_{ I^\a \ph^j \ph^l }=0,
\label{2.13}\\
&&  \Ga_{I^\a \ph^j}^{\ph^l}, \;
    \Ga_{I^\a I^\b}^{I^\ga} 
    \; \mbox{{\em{(}}and thus also } \Ga_{I^\a I^\b \ph^j}
    \mbox{{\em{)}} do not depend on } \ph^m,
\label{2.14} \\
&&  \Ga_{I^\a I^\b}^{\ph^j}, \;
    \; \mbox{{\em{(}}and thus also } \Ga_{I^\a I^\b I^\ga}
    \mbox{{\em{)}} are liner in } \ph^m.
\label{2.15} 
\eeqa
$\forall j, k, l, m, \a, \b,$ and $\ga = 1, \ldots, n.$
\eth

Here and later Greek letters will be used as indices for action 
coordinates and Latin letters as indices for angle coordinates. 
\hfill \\
\proof We construct a connection with the required properties 
by chart extension. Suppose that we have defined
the connection on a union $\pi^{-1}(\U)$ of action--angle charts 
$\pi^{-1}(U)$ \eqref{2.5} and we
would like to extend it to $\pi^{-1}(\ol{U}),$ with
action--angle coordinates as in \eqref{2.6}.
The Christoffel symbols for $\Ga$ on it will be defined consecutively in
the order listed in the statement of the theorem.

(1) Because of the transformation rules \eqref{2.7}--\eqref{2.8}
the symbols $\Ga_{\bp^j \bp^l }^{ \bI^\a}$ of $\Ga$ 
should satisfy
\[
\Ga_{\bp^j \bp^l }^{ \bI^\a} = 
   \sum_{j', l', \a'}
   \frac{\p \bI^\a}{\p I^{\a'}}
   \frac{\p \ph^{j'}}{\p \bp^j}
   \frac{\p \ph^{l'}}{\p \bp^l}
   \Ga_{ \ph^{j'} \ph^{l'} }^{ I^{\a'}}
\]
on $\pi^{-1}(U \cap \ol{U}),$  for $U \subset \U.$ But $\Ga$
satisfies \eqref{2.12} on $\pi^{-1}(U)$, so  we are free to put
$\Ga_{\bp^j \bp^l }^{ \bI^\a}=0$ on $\pi^{-1}(\ol{U}).$

(2) Similarly to part (1) we can define  
\[\Ga_{ \bp^j \bp^l }^{ \bp^m }=
   \Ga^{\bI^\b}_{\bI^\a \bp^j }=0,
\]
on $\pi^{-1}(\ol{U})$ and this will be consistent with the
definition of $\Ga$ on $\pi^{-1}(\U).$

(3) The transformation property that 
$\Ga_{\bI^\a \bp^j}^{\bp^l}$ should satisfy is
\beqa
\Ga_{\bI^\a \bp^j}^{\bp^l}
   &=& \sum_{\a' j' l'}
       \frac{\p I^{\a'}}{\p \bI^\a}
       \frac{\p \ph^{j'}}{\p \bp^j}
       \frac{\p \bp^{l}}{\p \ph^{l'}}
       \Ga_{I^{\a'} \ph^{j'}}^{\ph^{l'}} \nn    \\
   &+& \sum_{\a' j' \b'}      
       \frac{\p I^{\a'}}{\p \bI^\a}
       \frac{\p \ph^{j'}}{\p \bp^j}
       \frac{\p \bp^{l}}{\p I^{\b'}}
       \Ga_{I^{\a'} \ph^{j'}}^{I^{\b'}}  \nn     \\
   &+& \sum_{\a' j' \b'}
       \frac{\p \ph^{m'}}{\p \bI^\a}
       \frac{\p \ph^{j'}}{\p \bp^j}
       \frac{\p \bp^{l}}{\p \ph^{l'}}
       \Ga_{\ph^{m'} \ph^{j'}}^{\ph^{l'}} \nn    \\
   &+& \sum_{\b' j' m'}
       \frac{\p \ph^{m'}}{\p \bI^\a}
       \frac{\p \ph^{j'}}{\p \bp^j}
       \frac{\p \bp^{l}}{\p  I^{\b'}}
       \Ga_{\ph^{m'} \ph^{j'}}^{I^{\b'}}  \nn    \\
   &+& \sum_{l'} 
       \frac{\p^2 \ph^{l'}}{\p \bI^\a \p \bp^j}
       \frac{\p \bp^l}{\p \ph^{l'}}
\label{2.16}
\eeqa
on $\pi^{-1}(U \cap \ol{U}),$  for $U \subset \U.$
Eqs. \eqref{2.12}--\eqref{2.14} imply that 
on $\pi^{-1}(U \cap \ol{U})$ 
  $\Ga_{\ph^{m'} \ph^{j'}}^{I^{\b'}}=$      
  $\Ga_{I^{\a'} \ph^{j'}}^{I^{\b'}}=$
  $\Ga_{\ph^{m'} \ph^{j'}}^{\ph^{l'}}=0$       
and $\Ga_{I^{\a'} \ph^{j'}}^{\ph^{l'}}$ do not depend on     
$\ph.$ From this it is clear that the RHS of \eqref{2.16}
is a function on $\pi^{-1}(U \cap \ol{U}),$ which does not depend on
$\ph$ (see also \eqref{2.7}--\eqref{2.8}). It can be used to define
$\Ga_{\bI^\a \bp^j}^{\bp^l}$ on $\pi^{-1}(\U \cap \ol{U}),$  and
at the end we can define them on $\pi^{-1}(\ol{U})$ by any
continuation from the closure of $\pi^{-1}(U \cap \ol{U}).$ 

One deals with the symbols $\Ga_{\bI^\a \bI^\b}^{\bI^\ga}$ in the same
way.

(4) To define $\Ga_{\bI^\a \bI^\b}^{\bp^j}$ on
$\pi^{-1}(\ol{U})$ so that \eqref{2.15} holds, one just repeates the
construction from (3).

 {} From this construction we can get a torsion free symplectic connection
satisfying \eqref{2.12}--\eqref{2.15}. This can be done in two ways.
First, it is easy to keep on each step the symbols $\Ga_{x y z}$
completely symmetric. The second option is to construct a connection
$\Ga$ with the stated properties for $\Ga_{x y}^{z}$ and then to
symmetrize
$\Ga_{x y z}.$ The result is a connection with all the properties from
\thref{2.1}.
\hfill \qed

\bre{2.2} One can interprete properties 
\eqref{2.12}--\eqref{2.14} in a more geometric language. 

For each
$x\in M, b \in B,$ $T_x M^b \subset T_x M$ can be canonically
identified
with $T^*_b B,$ by linearization of the action of $T^*_bB/P_b$ on
$M^b$. Denote this identification by
$$
\i_x: \; T_x M^b \ra T^*_b B.
$$
Let $\Ga_B$ and $\Ga_M$ be two connections on $B$ and $M$ 
respectively with the following property.

If $y \in M_{b'}, b' \in B$ and $\ga$ is a curve connecting
$x$ and $y$ then 
\beq
\t_\ga^M(e) = \i^{-1}_y
                 \left( \t_{\pi(\ga)}^B(\i_x(e))\right),
\label{2.tr}
\eeq
where  $\t^M$ and $\t^B$ denote the parallel translations of tangent
vectors to $M$ and cotangent vectors to the base $B$ along the specified
curves (corresponding to $\Ga_M$ and $\Ga_B).$

Such a connection $\Ga_M$ satisfies  properties \eqref{2.12}--\eqref{2.13} 
because for $\ga \subset M^b$ 
$$
\t_\ga^B \equiv \i_y^{-1} \i^x: T_x M^b \cong T_y M^b,
$$ 
i.e. $M^b$ is a flat submanifold. Property \eqref{2.14} also holds for
$\Ga_M,$ because the parallel translation $\t^M$ of vectors in
$T_x M^b$ comes from the translation $\t^B$ on the base
(see \eqref{2.tr}).

Eq. \eqref{2.tr} can be actually used to define $\Ga_M$ on
$T_x M^b,$ and it will satisfy \eqref{2.12}--\eqref{2.14}.
This argument can substitute parts (1)--(3) of the above proof. 
\ere

\bre{2.3} $P$ is a covering of $B$ and we can define parallel translation
of the vectors of the lattice $P_b \subset T^*_bB$ along any curve in $B.$
Because of the nature of $P$ $(P_b$ is the kernel of the action 
$\tau$ \eqref{2.2} of $T^*_bB$ on $M_b)$ this can be extended by linearity
to parallel translation of vectors in $\span P_b \subset T^*_bB.$ In
the special case when $\pi$ has compact fibers $(k=n)$ it gives a torsion
free flat connection on $B.$
Using it as described in the previous remark we get a connection on $M$
for which all Christoffel symbols from \eqref{2.12}--\eqref{2.14} vanish.
\ere

%%%%%%%%%%%%%%%%%%%%%%%%%%%%%%%%%%%%%%%%%%%%%%%%%%%%%%%%%%%%%%%%%%%%%%%%%%%%%%%
\sectionnew{Existence of star products preserving a Lagrangian 
fiber bundle}
%%%%%%%%%%%%%%%%%%%%%%%%%%%%%%%%%%%%%%%%%%%%%%%%%%%%%%%%%%%%%%%%%%%%%%%%%%%%%%%

The goal of this section is to prove the existence of nontrivial
deformations of $\CM,$ that keep $\pcb$ undeformed. In the next section we
show that up to equivalence of star products these are all deformations of
$\CM$ in which $\pcb$ generates a commutative subalgebra.

\bth{3.1} Let $\Om \in Z^2(M)[[\h]]$ be a formal deformation of the
symplectic form $\om:$
\beq
  \Om = \om + O(\h),
\label{3.1}
\eeq
such that $\pi: \: M \ra B$ is a Lagrangian fibration
with respect to $\Om,$ that is 
\beq
  \Om \Big|_{\pi^{-1}(b)}=0, \;
  \forall b \in B.
\label{3.2}
\eeq
Then there exists a Fedosov's star product with characteristic class
$[\Om]$ that keeps $\pcb$ undeformed.
\hfill
\eth

Denote by $\overline{W}_B$ the Weyl bundle on $B$ and by $\Ga(\WB)$
the pull back to $M$ of the space of its global sections:
\beq
\Ga(\WB) = \pi^*(\Ga(\overline{W}_B)) \subset \Ga(\WM)
\label{3.3}
\eeq
Consider also the conormal bundle to the fibration $\pi$
\[
(N^*_{\pi})_x = \{ v \in T_x^*M | <v, e>= 0,\: \forall e \in T_x
M^b\} \quad
\mbox{for } x \in M^b, \; b \in B.
\]
We will call 
\beq
\Wp = S(N^*_\pi) \subset \WM
\label{3.4}
\eeq
Weyl bundle of the fibration $\pi.$ 
Also denote $\WBh =\WB[[\h]],$ $\Wph=\Wp[[\h]],$
and $\WMh=\WM[[\h]].$ The following inclusions are obvious 
\beq
\Ga(\WB) \subset \Ga(\Wp) \subset \Ga(\WM)
\label{3.5}
\eeq
and clearly $\Ga(\Wp)$ is an undeformed subalgebra of $\Ga(\WMh)$ with
respect to the fiberwise Moyal--Weyl product \eqref{1.3}.

Locally, in action--angle coordinates $(I, \ph)$ on some $\pi^{-1}(U)$ as
in eq. \eqref{2.5} denote by $J^\b$ and $\psi^j$ the images of $ d I^\b$
and $\ph^j$ in $\WM = S(T^*M)$ for $\b, j = 1, \ldots, n.$ In this
notation $\Ga (\Wph, \pi^{-1}(U))$ consists of sections \eqref{1.loc}
that are power series in $J^\b$ only and 
$\pi^*(\Ga (\overline{W}_{B, \h}, U))$
of the
subspace of sections for which in addition
$a_{k, \a}$ depend on the action coordinates $I^\b$ only.

In Subsect. 4.1 we show that for any closed two-form $\Om$ as in the
statement of \thref{3.1} there exists a closed form $\Om'$ of a special
type cohomologous to it. In the second subsection we deduce \thref{3.1}
  {} from a stronger result stating that Fedosov's construction 
applied to $\Om'$ exponentiates $\pcb$ to a subalgebra of $\Ga(\WBh).$
The latter, as we remarked above, is an undeformed subalgebra of 
$(\Ga(\WMh), \circ).$

\subsection{Construction of special closed two-forms}
Denote by $\Ga(\WMh\otimes\La)^\pol$ the subalgebra of
$\Ga(\WMh \otimes \La)$ consisting of those sections which in any
action--angle chart \eqref{2.5} are formal power series
\beq
\sum_{l=0}^{\infty}\sum_{\a, \ga}
\h^l f_{a, \ga, l}(I, \ph)y^\a d x^\ga
\label{3.6}
\eeq
with $f_{a, \ga, l}(I, \ph)$ being polynomials in the variables
$\ph^{k+1}, \ldots, \ph^{n},$ where $x=(I, \ph),$ $y=(J, \psi),$ 
and $d x^\ga= d x^{\ga_1} \wedge \ldots \wedge d x^{\ga_s}$
(see also the remarks to eq. \eqref{2.5}). 
The above requirement does not depend on the choice of a particular 
action--angle coordinate system $(I, \ph)$ on an open $\pi^{-1}(U)$
because of the transformation rules \eqref{2.7}--\eqref{2.8}.
The same formulas make possible to define a filtration on
$\Ga(\WMh \otimes \La)^\pol$ by letting (locally)
\beq
  \deg \ph^j = 
  \deg d \ph^j =
  \deg \psi^j = 1,
  \quad
  j = 1, \ldots, n.
\label{3.7}
\eeq
Denote the part of $\Ga(\WMh \otimes \La)^\pol$ of $\deg \leq l$ by
\beq
\Ga(\WMh \otimes \La)^{(l)}
\label{3.8}
\eeq
(not to be mistaken with the grading \eqref{1.grad}).
In the same way we define the subspaces 
$\Ga(\WMh)^\pol$ and $\Ga(\La)^\pol$ of
$\Ga(\WMh \otimes \La)^\pol$ and restrict
the filtration \eqref{3.8} to them.

We are ready to state the main result of this subsection.

\bpr{3.2} Let $\Om \in Z^2(M)$ and 
\[
\Om \Big|_{\pi^{-1}(b)}=0
\]
for all $b \in B.$ Then there exists $\Om' \in Z^2(M)$ cohomologous
to $\Om$ such that
\beq
\Om' \in \Ga(\La^2)^{(1)}.
\label{3.10}
\eeq
{\em{(}}Then also automatically $\Om' |_{\pi^{-1}(b)}=0, \;
\forall b \in B.${\em{)}}
\epr

Both the above result and the Lemma that we need for its proof are
well known but we will sketch two elementary proofs of them for 
completeness.

\ble{3.3}{\em{(}}i{\em{\/)}} Let $\a \in Z^2(U \times \Rset^n)$ for an
open domain $U$ in $\Rset^n.$ Assume that
\beq
  \a \Big|_{\pi_1^{-1}(u)}=0, \; \forall u \in U,
\label{3.11}
\eeq
where $\pi_1: \: U \times \Rset^n \ra U$ is the projection on the first
factor. Then there exists a one-form $\b$ such that
\beq
  \b \Big|_{\pi_1^{-1}(u)}=0, \; \forall u \in U
  \quad \mbox{and} \quad d \b = \a.
\label{3.12}
\eeq

{\em{(}}ii{\em{\/)}} If $\b \in Z^1(U \times \Rset^n)$ and 
\[
  \b \Big|_{\pi_1^{-1}(u)}=0, \; \forall u \in U,
\]
for an open subset $U \subset \Rset^n$ then 
\[
\b = \pi^*(\b_1), 
\]
for some $\b_1 \in Z^1(U).$
\ele
\proof Part (i) Denote by $(I^1, \ldots, I^n)$ and $(\ph^1, \ldots,
\ph^n)$ the coordinates on $U$ and $\Rset^n$ respectively. If
\[
   \a = f_{j l}(I, \ph) d I^j\wedge d I^l 
      + g_{j l}(I, \ph) d \ph^j\wedge d I^l
\]
then $\a^0 =f_{j l}(I, 0) d I^j\wedge d I^l$ is closed in $U$ and there
exists a one-form $\b^0 = a_l(I) d I^l$ such that $d \b^0 = \a^0.$
Using the
assumption that $\a$ is closed it is easy to prove that
\[
  \b = \b^0
     + \left( \int_0^\ph
       g_{jl}(I, \phi) d \phi^j
       \right) 
       d I^l
\]
satisfies \eqref{3.12}. 

Part (ii) If $\b= b_l(I, \ph) d I^l$ then $d \b = 0$ implies 
$d_\ph b_l(I, \ph)=0,$ $\forall I \in U.$ Thus $b_l(I, \ph)$ depends on
$I$ only and we can put $\b_1 = b_l(I) d I^l.$
\hfill
\qed \\
\hfill \\
{\em{Proof of \prref{3.2}.}} We construct a two-form $\Om',$ satisfying
\eqref{3.10} and a one-form $\ga$ such that $\Om -\Om'= d \ga$
by chart extension.

Let $(I, \ph)$ be an action--angle chart on some 
$\pi^{-1}(U) \cong U \times(\Rset^n/\Zset^k).$
Denote 
\[
\Om_U= \Om \Big|_{\pi^{-1}(U)}
\]
and let its pull back to $U \times \Rset^n$ be
\beq
   \wt{\Om}_U = f_{j l}(I, \ph) d I^j\wedge d I^l
      + g_{j l}(I, \ph) d \ph^j\wedge d I^l
\label{3.15}
\eeq
(here $f_{jl}$ and $g_{jl}$ are $\Zset^k$ periodic in 
$\ph^1, \ldots, \ph^k).$
\leref{3.3} applied to $\wt{\Om}_U$ gives a one-form
$\tilde{\b}_U = a_l(I, \ph) d I^l$ such that 
$d \tilde{\b}_U = \wt{\Om}_U.$
For any $l$ and $I=c=(c^1, \ldots, c^n)= \const$
\[
g_{jl}(c, \ph) d \ph^j \in Z^2(\Rset^n/\Zset^k)
\]
because $\Om$ is closed. $H^2(\Rset^n/\Zset^k)\cong \Rset^k$ and there
exists a one-form with coefficients that do not depend on the angle
variables $\ph$
\beq
   \sum_j  a'_{jl}(c) d \ph^j \sim g_{j l}(c, \ph) d \ph^j 
\label{3.19}
\eeq
where the first sum is over $1 \leq j \leq k,$ i.e. just coordinates 
$\ph^1, \ldots, \ph^k$ related to the lattice $\Zset^k$ appear. Then
\beq
  \ga_U := (a_l(I, \ph) -\sum_j a'_{jl}(I)\ph^j) d I^l
\label{3.20}
\eeq
descends to a one-form  on $U \times \Rset^n/\Zset^k,$ and so does 
\beq
  \Om'_U := \Om_U - d \ga_U = d( (\sum_j a'_{jl}(I) \ph^j) d I^l) 
\label{3.21}
\eeq
(to a two-form on $U \times \Rset^n/\Zset^k).$
Clearly
\[
  {\Om'}_U \in \Ga(\La^2, \pi^{-1}(U) )^{(1)}.
\]

Suppose now that we have already constructed $\Om'$ and $\ga$ on a union
of action--angle charts and we want to extend them to an open
$\pi^{-1}(U)$ on $M$ as in \eqref{2.5}. The intersection
of $\pi^{-1}(U)$ with this union is $\pi^{-1}(U^0),$ for some 
open $U^0 \subset U.$ Denote the restriction of $\ga$ to 
$\pi^{-1}(U^0) \cong U^0 \times \Rset^n/ \Zset^k$ by
$\ga_{U^0}.$ The second part of \leref{3.3} 
together with the transformation rules \eqref{2.7}--\eqref{2.8}
imply that
\beq
  \ga_{U^0} = (a_l(I, \ph) - a'_{jl}(I)\ph^j -d_l(I)) d I^l 
\label{3.ga}
\eeq
for some functions  $d_l(I)$ on $U^0$ (the sum over $j$ is from 1 to
$k$ again). 
%\\ $(b_{1l}(I), \ldots, \b_{kl}(I))$ represents the
%cohomology class of $d_\ph a_l(I, \ph)$ under the identification
%$H^2(\Rset^n/\Zset^k) \cong \Rset^k.$ The functions $a_{jl}(I)$ from eq.
%\eqref{3.20} have the same property (see \eqref{3.19}), so
%$b_{jl}(I)= a_{jl}(I)$ for $I \in U^0.$ 
Finally we can extend $d_l(I)$ to $C^\infty$-functions on $U$ from the
closure of $U^0$ and define an extension $\ga_U$ of $\ga_{U^0}$ to 
$\pi^{-1}(U)$ by the same formula \eqref{3.ga}.
${\Om'}_{U^0}$ also can be extended to $\pi^{-1}(U)$ by
${\Om'}_U =\Om_U - d \ga_U.$ Locally on $\pi^{-1}(U)$
\[
  \Om'|_U= d(a_{jl}(I) \ph^j +d_l(I))\wedge d I^l
  \in \Ga(\La^2, \pi^{-1}(U))^{(1)}.
\]
\qed

\subsection{Construction of star products preserving a Lagrangian fiber
bundle}
The next theorem is the main result in this subsection.
It makes \thref{3.1} more precise, the latter is a consequence
of it due to \prref{3.2}.

\bth{3.4} Let $\Om' \in Z^2(M)[[\h]]$ be a formal deformation of $\om,$
$(\Om' = \om + O(\h))$ with the property from \prref{3.2}
\beq
\Om' \in \Ga(\La^2 M[[\h]])^{(1)}.
\label{3.gr}
\eeq
Then any Fedosov's connection $D$ with curvature $\Om'$ based on a 
symplectic connection $\p$ with the properties of \thref{2.1} exponentiates
$\pcb[[\h]]$ to a subalgebra of $\Ga(\WBh):$
\[
  \sig_D \left( \pcb[[\h]]
         \right) 
  \subset
  \Ga(\WBh).
\]
\eth
We will need a simple Lemma.

\ble{3.5}Assume that $a \in \Ga(\WMh)^{(q)}$ and
$b \in \Ga(\WMh)^{(r)},$ then

{\em{(}}i{\em{\/)}} $[a, b] \in \Ga(\WMh)^{(q+r -1)},$
 
{\em{(}}ii{\em{\/)}} $\p a \in \Ga(\WMh)^{(q)},$ for any connection $\p$
with the
properties of \thref{2.1},

{\em{(}}iii{\em{\/)}} $\de a$ and $\de^{-1} a \; \in \Ga(\WMh)^{(q)}.$
\ele
\proof Parts (i) and (iii) are trivial. Part (ii) follows from
\beqa
   \p I^\a  &=& \Ga_{I^\b I^\ga}^{I^\a} J^\b d I^\ga \;
          \in \Ga(\WMh)^{(0)}, \nn \\        
   \p \ph^l &=& \Ga_{I^\b \ph^j}^{\ph^l} J^\b   d \ph^l 
           + \Ga_{\ph^j I^\b}^{\ph^l} \psi^j d I^\b 
           + \Ga_{I^\a I^\b}^{\ph^l}  J^\a   d I^\b \;
           \in \Ga(\WMh)^{(1)}
\nn
\eeqa
for any connection $\p$ satisfying \eqref{2.12}--\eqref{2.15}.
\hfill
\qed   \\
{\em{Proof of \thref{3.4}.}} Eq. \eqref{1.11} implies that the curvature
form $R$ of $\p$ (see \eqref{1.15}) in local Darboux coordinates is
\[
  R = \p \Ga +\frac{1}{\h}(\Ga \circ \Ga),
\]
where $\Ga= (1/2)\Ga_{ijl} y^i y^k d x^l$ is the connection
form of $\p$ in the considered coordinate system. Clearly from
\eqref{2.12}--\eqref{2.15} we have
\[
  \Ga \in \Ga((W\otimes\La^1)^{\pol}, \pi^{-1}(U))^{(1)}
\]
and from \leref{3.5} we get that 
\[
  R \in \Ga(\WM\otimes\La^2)^{(1)}
\]
(see \leref{3.5}). Using the same Lemma again from eq. \eqref{1.18} we 
obtain
\beq
\ga \in \Ga(\WMh\otimes\La^1)^{(1)}.
\label{3.ga1}
\eeq
Finally \leref{3.5} applied to \eqref{1.20} gives  
\[
  \sig(f) \in
  \Ga(\WMh)^{(0)} = \Ga(\WBh), \; \forall f \in \pcb.
\]
\qed

\subsection{Another proof of the existence result}
Here we sketch another proof of \thref{3.1} that does not use
\prref{3.2}. The idea is to use Fedosov's construction with 
$\Om= \om + O(\h) \in H^2(M)[[\h]]$ instead of $\om.$

Consider $M$ as a $C^\infty[[\h]]$ manifold. That is we take 
$\O^\h_M=\O_M[[\h]]$ as structure sheaf on $M$ ($\O_M$ is the sheaf
of $C^\infty$ functions on $M)$. For any open $U \subset M$
$\Ga(\O^h_M, U)=\Ga(\CM, U)[[\h]]$ and a coordinate system on $U$ is 
specified by the images of the coordinate functions $x_1, \ldots, x_{2n}$
under an isomorphism
\[
F:
C^\infty(U_o)[[\h]] \ra \Ga(\O^h_M, U),
\]
for an open $U_o \subset \Rset^{2n}.$ Any composition
$f(g_0 + \h g_1 + \dots)$ for some $C^\infty$ functions $f, g_0, g_1,
\ldots$ is computed using the Taylor rule 
$f(g_0 + \h g_1 + \dots)= \sum_{l=0}^\infty f^{(l)}(g_0)G^l/l!,$
$G= \h g_1 + \h^2 g_2 + \dots.$ One can consider $(M, \Om)$ as a
symplectic $C^\infty[[\h]]$ manifold and prove the existence of a
covering of it with action--angle charts. Let $\pi: M \ra B$ be a
fiber bundle which is Lagrangian  with respect to $\Om$ in the usual
sense. Using \leref{3.3} one constructs a symplectic connection
$\ol{\Ga}$ on $(M, \Om)$ with covariant derivative $\bar{\p}: \Ga(TM)
\ra
\Ga(TM)[[\h]],$ that satisfies properties \eqref{2.12}--\eqref{2.15} for
the data $M,$ $\pi,$ $\Om.$
It is obvious that the constant term of $\bar{\p}$ 
\[
\p: \; TM \ra TM
\]
is a connection with the properties of \thref{2.1} for the data
$M,$ $\pi,$ $\om.$ 

Now we can run Fedosov's construction with $\om$ and $\p$ replaced by
$\Om$ and $\bar{\p}.$ $\WMh$ will be equiped with the Moyal--Weyl
product corresponding to $\Om$, let us denote it by $\bar{\circ}.$
Consider a Fedosov's connection for $(\WMh, \bar{\circ})$
\[ 
\ol{D}= \bar{\p} - \de + \frac{1}{\h}\ad(\bar{\ga}).
\] 
with curvature form (same) $\Om$. As in the proof of \thref{3.4}
(see also \leref{3.5}) one proves that it induces a star product $\st$ on
$M$ which keeps $\pcb$ undeformed. In the remaining part of this Remark we
show that $\st$ has characteristic class $[\Om].$ 

One can find a fiberwise linear map $F: TM \ra TM[[\h]],$
such that $F(\om)=\Om.$ Then
\[
F^*: (\WMh, \bar{\circ}) \ra (\WMh, \circ)
\]
is an isomorphism. Continue $F^*$ to $\WMh\otimes\La$ 
trivialy on the second factor. Clearly 
\beqa
&& F^* \bar{\p} (F^*)^{-1}= \p + \ad_\circ(\a) \;  \mbox{and} 
\nn \\
&& F^* \de (F^*)^{-1}= \de + \ad_\circ(\b)
\nn 
\eeqa 
for some one-forms $\a$ and $\b$ with values
in $S^2(T^*M)$ and $S^1(T^*M)$ respectively. 
Here $\ad_{\bullet}$ refers to the the Lie bracket associated with the
product $\bullet.$ As a consequence of this we get that 
\beqa
D &:=& (F^*)^{-1} \ol{D} F^*
\nn \\
  &=& \p-\de + \frac{1}{\h}\ad_\circ(F^*(\bar{\ga}) + \h \a 
                                                     - \h \b).
\nn
\eeqa
is a Fedosov's connection for $(\WMh, \circ).$ By straightforward
computation one shows that the Weyl curvature of $D$ is $\Om$ and
$F^*$ induces isomorphism between $W_{\ol{D}}$ and $W_{D}.$
It is now clear that the star products coming from $\ol{D}$
and $D$ are equivalent and have characteristic class $[\Om].$
\hfill 
%%%%%%%%%%%%%%%%%%%%%%%%%%%%%%%%%%%%%%%%%%%%%%%%%%%%%%%%%%%%%%%%%%%%%%%%%%%%%%%
\sectionnew{Classification of star products preserving a Lagrangian
fiber bundle}
%%%%%%%%%%%%%%%%%%%%%%%%%%%%%%%%%%%%%%%%%%%%%%%%%%%%%%%%%%%%%%%%%%%%%%%%%%%%%%%
The results of this Section are summarized in the following theorem.

\bth{4.1} Let $\st$ be a star product on the symplectic manifold $(M,
\om)$ for which $\pcb$ generates a commutative subalgebra of $\CM[[\h]].$
Then its characteristic class has a representative $\Om \in Z^2(M)[[\h]]$
such that $\pi: M \ra B$ is Lagrangian with respect to $\Om:$
\beq
  \Om \Big|_{\pi^{-1}(b)}=0,
  \quad \forall b \in B.
\label{4.1}
\eeq
In addition, there is always a Fedosov's star product on $(M, \om)$ that
is
equivalent to $\st$ and keeps $\pcb$ undeformed.
\eth

The first step towards our proof of \thref{4.1} is a version of a well
known fact, which states that any commutative deformation of the algebra
of functions on a smooth manifold is equivalent to a trivial
one.

\ble{4.2}Let $\st$ be a star product on $(M, \om)$ for which the algebra
generated by $\pcb$ is commutative. Then there exists an equivalent
product $\pe$ on $(M, \om)$ that keeps $\pcb$ undeformed.
\ele
We will skip the proof of this lemma as it is similar to the proof of
\prref{4.4} that follows next.

The second step in our proof of \thref{4.1} is to exponentiate $\pe$
to a fiberwise star product $\pew$ on $\Wh$ using Emmrich--Weinstein
construction \cite{EW} reviewed in Subsect.~2.5.1 with
a symplectic connection $\p$ having the properties of \thref{2.1}.

\ble{4.3}Let $\pe$ be a star product on $(M, \om)$ that keeps $\pcb$
undeformed. Denote the fiberwise star product $\pew$ on $\WMh$
obtained from $\pe$ by the exponentiation construction of Subsect.~2.5.1
using a symplectic connection $\p$ with the properties of \thref{2.1}
{\em{(}}see eq. \eqref{1.Wprod}{\em{\/).}} Then $\pew$ keeps each fiber
of $\Wph$
undeformed.
\ele
\proof In the exponentiation procedure from Subsect.~2.5.1 one first
computes the one-form $\ga$ for the nonlinear flat connection $D^0$ 
on $W$ (see \eqref{1.22}). It is done inductively using \eqref{1.23}.
Similarly to the proof of \thref{3.4} (see especially \eqref{3.ga1})
one gets 
\[ 
  \ga \in \Ga(\WM)^{(1)}.
\]
Thus
\[
  \exp(f) \in \Ga(\WM)^{(0)}= \Ga(\WB)
  \quad \forall f \in \pcb
\]
(see eq. \eqref{1.25} for the definition of $\exp).$ The closure of
\beq
\span \{ \exp(f)_x | f \in \pcb\}
\label{4.3}
\eeq
in $S(T^*_xM)$ is $S((N^*_\pi)_x)=(\Wp)_x,$ $\forall x \in M,$ which
proves 
the Lemma. 
%Let us remark that in an action--angle chart $(I, \ph)$
%on some open $\pi^{-1}(U)$
%\beqa
%  &&S((N^*_\pi)_M) = \Rset[ J^1, \ldots, J^n],
%\label{4.4} \\
%  &&W_m = \Rset[J^1, \ldots, J^n, \psi^1, \ldots, \psi^n]
%\label{4.5}
%\eeqa
\hfill
\qed

The next step is to study the equivalence between the product $\pew$ 
and the Moyal--Weyl product $\circ$ \eqref{1.3}.

\bpr{4.4} There exists a leafwise differential operator $P$ on $TM$
of the type
\beq
P = \id + \h P_1 + \h^2 P_2 + \dots,
\label{4.6}
\eeq
for some differential operators $P_j$ on $TM$ with coefficients in $\WM$
such that 
\beq
P(a \pew b) = P(a) \circ P(b), \quad
\forall a, \: b \in (W)_x, \; x \in M
\label{4.7}
\eeq
and
\beq
P a = a, \quad 
\forall a \in (\Wp)_x \subset (W)_x.
\label{4.8}
\eeq
\epr
\proof Denote by $C^l$ the bundle over $M$ whose fiber at $x \in M$
consists of (local) Hochschild $l$-cochains for $(\W)_x=\J_0(T_x M)$
(equipped with the standard fiberwise commutative product). 
Any such cochain is represented by an $l$-differential operator on the
tangent space $T_x M.$

Let us define
\beqa
(C^l)^s = \{ B \in C^l |
&&B \: \mbox{is represented by an $l$-differential operator}
  \nn \\
&&\mbox{of degree $> 1$ in at lest one component or it is}
  \nn \\
&&\mbox{of degree 1 in any component and } A(B)=0 \},
  \label{4.9}
\eeqa
where $A$ denotes the complete antisymmetrization and
\beq
(C^l)^a = \{ B \in C^l | 
          B \: \mbox{is antisymmetric of degree 1 in each component}
          \}.
\label{4.10}
\eeq
Let $(C^l)^a_\pi,$ $(C^l)^s_\pi,$ and $C^l_\pi$
be the subbundles of $(C^l)^a,$ $(C^l)^s,$ and $C^l$
consisting of those cochains that vanish if one of their
arguments is in $\Wp.$ It is clear that
\beq
C^l_\pi = (C^l)^a_\pi \oplus (C^l)^s_\pi.
\label{4.11}
\eeq
By $\de$ and $[.,.]$ (in this proof only) we denote the
(fiberwise) differential and the (super) Lie bracket in the Hochschild
complex $C.$ $\a \in \Ga(C^2)$ will denote the cochain corresponding to
the fiberwise Poisson bracket $\{.,.\}_\fib$ on $\WM.$ The following
statements are well known:

(i) $\de((C^l)^s) \subset (C^{l+1})^s$ and
\[
\dots \stackrel{\de}{\ra}
      (C^{l-1})^s   \stackrel{\de}{\ra}
      (C^l)^s       \stackrel{\de}{\ra}
      (C^{l+1})^s   \stackrel{\de}{\ra}
\dots
\]
is acyclic. In particular
\[
      (C^1)^s \stackrel{\de}{\ra}
      (C^2)^s 
\]
is an inclusion onto $\ker \: \de_2((C^2)^s \ra (C^3)^s).$

(ii) $(TM, \{., .\})$ is a regular Poisson manifold with symplectic leaves 
$T_x M,$ $x \in M.$ The restriction of $\a$ to each leaf canonically
identifies the fiber of $(C^l)^s$ at $x$ with the space of forms
$\Om(T_xM)$ on $T_x M$ with coefficients in $(W)_x= \J_0 (T_xM).$
Under this identification $\{\a, .\}$ goes to the fiberwise
differential $d_{TM}$ (see \cite{NT}).

We prove that there exists a differential operator $P$ 
as in \eqref{4.6} satisfying \eqref{4.7} such that
\beq
P_n \in \Ga((C^2)_\pi).
\label{4.12}
\eeq
The last equation obviously implies that $P$ also satisfies \eqref{4.8}.
According to \eqref{4.11} $P_n$ should decompose as
\[
P_n = P_n^s + P_n^a,
\; \mbox{where} P_n^s \in   \Ga((C^1)^s_\pi).
\; \mbox{and}   P_n^a \in \Ga((C^1)^a_\pi)
\]
The operator $P$ is constructed iteratively. At the $n$-th step
assuming that \eqref{4.7} holds modulo $\h^{n-1},$ we define
$P_n^a$ and $P_{n-1}^s$ in such a way that \eqref{4.7} is true modulo
$\h^n.$ The equations for $P_n^s$ and $P_{n-1}^a$ are
\beqa
  &&\de P_n^s = Q_n^s, 
\label{4.14} \\
  &&[\a, P_{n-1}^a] = Q_{n-1}^a.
\label{4.15} 
\eeqa
with some $Q_n^s \in \Ga((C^2)^s_\pi)$ and 
$Q_{n-1}^a \in \Ga((C^2)^a_\pi),$
which are known from the previous steps and satisfy
\beqa
  &&\de Q_n^a = 0,
\label{4.16} \\
  &&[\a, Q_{n-1}^s] = 0.
\label{4.17}
\eeqa

According to (i) \eqref{4.14} has a unique solution $P_n^s$ in
$\Ga((C^1)^s)$ and this solutions is in $\Ga((C^1)^s_\pi).$ The existence
of a solution $P_{n-1}^a \in \Ga((C^1)^a_\pi)$ of \eqref{4.15} can
be proved in the following way. The identification from (ii) associates
with $Q_{n-1}^a$ a family of two-forms $\b^2_x$ on $T_x M, x \in M$
with coefficients in $\J_0(T_x M).$ The restriction of $\b^2_x$ to 
$T_x M^b$ vanishes for any $x \in M$ (here $M^b$ denotes the fiber of 
$\pi$ passing trough $x$). A formal analog of \leref{3.3} implies that
there exists a family of one forms $\b^1_x$ on $TM$ whose resetrictions
to $T_x M^b$ vanish and on each tangent space $T_x M$
\[
d \b^1_x = \b^2_x.
\]
Applying the identification from (ii) back to $\b^1_x$ plus some standard 
gluing arguments give the existence of a
solution to 
\eqref{4.15} of the desired form.
\hfill 
\qed

Finally \thref{4.1} follows directly from the next Theorem, which
makes it more precise.

\bth{4.6}Let $\st$ be a star product on $(M, \om)$ for which
$\pcb$ generates a commutative subalgebra in $\CM[[\h]].$
Then there exists an injective homomorphism
\beq
\ro: (\CM[[\h]], \st) \hra (\Ga(\WMh), \circ)
\label{4.18}
\eeq
such that 
\beq
\ro (\pcb) \subset \Ga(\WBh).
\label{4.19}
\eeq
There exists a Fedosov's connection $D$ on $\Wh$ for which  
\beq
\ro = \sig_D.\
\label{5.con}
\eeq 
Its Weyl curvature satisfies \eqref{4.1}. The corresponding
Fedosov's product is equivalent to $\st$ and keeps $\pcb$ undeformed.
\hfill
\eth
\proof Most of the work needed is already done. The homomorphism
$\ro$ is just the composition of the maps from Lemmas
\ref{l4.2}, \ref{l4.3}, and \prref{4.4}. 
It has properties (i) and (ii) of Subsect.~2.5.1
and from \thref{1.5} we get that there exists a Fedosov connection 
$D= \p -\de +(1/\h)\ad(\ga)$ satisfying \eqref{5.con}.
Then for any $f \in \CM$
\[
\left(
  \p - \de - \frac{1}{\h} \ad(\ga)
\right)
\ro(f) = 0, \quad
\forall f \in \CM.
\]
We have the freedom to choose the initial connection $\p.$ If it has
the properties of \thref{2.1}, then
\[
(\p -\de) \ro(f) \in \Ga(\WMh\otimes\La^1)^{(0)},
\quad \forall f \in \pcb, 
\]
and thus
\[
[\ga, \ro(f)] \in \Ga(\WMh\otimes\La^1)^{(0)},
\quad \forall f \in \pcb.
\]
The operator $P$ from \prref{4.4} acts
trivially on $\pcb,$ so
\[
\ro(f)=\exp(f) \in \Ga(\Wph), 
\quad \forall f \in \pcb,
\]
where $\exp$ is the exponential map from the proof of \leref{4.3}. 
This proves the fact about Fedosov's product in \thref{4.6}. It also
implies that the completion of
\[
\span \{ \ro(f)_x |
             f \in \pcb[[\h]] \}
\]
in $(\Wh)_x$ is $(\Wph)_x$ $(\forall x \in M)$ and therefore that for
any open domain $U \subset B$
\beq
\ga\big|_U \in \Ga(\WMh\otimes\La^1, \pi^{-1}(U))^{(1)} +
           \Ga(\Wph\otimes\La^1, \pi^{-1}(U) ).
\label{4.deco}
\eeq
It is clear that in local action--angle coordinates $(I, \ph)$
\[
 (\Wp)_m = \Rset [[J^1, \ldots, J^n]]
\]
and that it is the commutant of itself in
\[
\left( (\WMh)_m= \Rset[[J^1, \ldots, J^n, \psi^1, \ldots, \psi^n]], \circ
\right).
\]
Recall that by $J^l$ and $\psi^l$ we denote the images of $d I^l$ and
$d \ph^l$ in $(W)_x$ respectively, $x \in M.$ At the end the Weyl
curvature of $D$ is \beqa
\Om &=& R + \p \ga - \de \ga +\frac{1}{\h}\ga^2
   \nn \\
    &=&(\ga \circ \ga)_0
   \nn \\
\eeqa
(see the discussion before \thref{1.2} for the definition of constant term
$(.)_0.)$ Locally write $\ga$ as a sum of two terms according to
\eqref{4.deco}.
Then nonzero contribution to $(\ga \circ \ga)_0$ give terms of the
following type only:
\beqa
&&\ga'_1 \circ \ga''_1, \;
\mbox{for } \ga'_1, \: \ga''_1 \in \Ga(\WMh^j\otimes\La^1)^{(1)},
j = 1 \: \mbox{or} \: 2,
\label{4.c1} \\
&& \ga_1 \circ \ga_2,
\; \mbox{for } \ga_1 \in \Ga(\WMh^1\otimes\La^1)^{(1)},
 \: \ga_2 \in \Ga(\Wph\otimes\La^1).
\label{4.c2}
\eeqa
In case \eqref{4.c1}
\[
\ga'_1 \circ \ga''_1 \in \Ga(\La^2)^{(2)}[[\h]]
\]
and all forms in $\Ga(\La^2)^{(2)}[[\h]]$
satisfy \eqref{4.1}.

Consider case \eqref{4.c2}. Locally $\ga_1$ is a polynomial in $\psi^j$ of
degree at most 1. The product \eqref{4.c2} eats out the $\psi$'s and what
remains is in $\Ga(\La, \pi^{-1}(U))^{(0)}[[\h]].$ As a one-form it
vanishes on $M^b, \: \forall b \in B$ and so does $\ga_1 \circ \ga_2.$
This ends the proof of \thref{4.6}.
\hfill \qed
%%%%%%%%%%%%%%%%%% References %%%%%%%%%%%%%%%%%%%%%%%%%%%%%%%%%%%%%%%%%%%%%%%
    
%%%%%%%%%%%%%%%%%%%%%%%%%%%%%%%%%%%%%%%%%%%%%%%%%%%%%%%%%%%%%%%%%%%%%%%%%%%%%%%
%%%%%%%%%%%%%%%%%%%%%%%%%%%%%%%%%%%%%%%%%%%%%%%%%%%%%%%%%%%%%%%%%%%%%%%%%%%%%%
\end{document}